\documentclass[12pt]{amsart}

\usepackage{amsmath,amssymb,amsthm}
\usepackage[dvipdfm]{graphicx}
\usepackage[dvips]{color}

\newtheorem{Thm}{Theorem}[section]
\newtheorem{Lem}[Thm]{Lemma}
\newtheorem{Cor}[Thm]{Corollary}
\newtheorem{Prop}[Thm]{Proposition}

\newtheorem{Claim}[Thm]{Claim}
\newtheorem{Property}[Thm]{Property}

\theoremstyle{remark}

\theoremstyle{definition}
\newtheorem{Def}[Thm]{Definition}
\newtheorem{Step}{Step}

\newcommand{\Shpec}{\mathop{\mathit{Spec}}\nolimits}

\newcommand{\Supp}{\mathop{\mathrm{Supp}}\nolimits}

\newcommand{\DF}{\mathop{\mathrm{DF}}\nolimits}

\newcommand{\NKLT}{\mathop{\mathrm{NKLT}}\nolimits}
\newcommand{\NLC}{\mathop{\mathrm{NLC}}\nolimits}
\newcommand{\NSLC}{\mathop{\mathrm{NSLC}}\nolimits}

\newcommand{\NN}{\mathop{\mathrm{NN}}\nolimits}

\newcommand{\cond}{\mathop{\mathrm{cond}}\nolimits}

\newcommand{\GL}{\mathop{\mathrm{GL}}\nolimits}

\begin{document}

\title[Stability via discrepancy]
{The GIT stability of Polarized Varieties via Discrepancy} 

\author{Yuji Odaka}
\dedicatory{Dedicated to Professor Shigefumi Mori on his Kanreki (60th birthday)} 
\address{Research Institute for Mathematical Sciences (RIMS), 
Kyoto University, Oiwake-cho, Kitashirakawa, Sakyo-ku, Kyoto
606-8502, Japan}
\address{{\textit{Current address}}:  
Imperial College London, 
South Kensington Campus, 
London SW7 2AZ, 
England}
\email{y.odaka@imperial.ac.uk}


\begin{abstract}
We prove that {\color{black}{various}} GIT semistabilit{\color{black}{ies}} 
of polarized varieties impl{\color{black}{y}} semi-log-canonicity. 
\end{abstract}

\maketitle


\section{Introduction}

For the study of the moduli of polarized varieties, 
Geometric Invariant Theory \cite{Mum65} 
(GIT, for short) is an important basis, 
because it constructs the moduli spaces as quotient schemes of the Hilbert schemes. 
In that theory, we must put restrictions on the objects to classify, which we call \textit{stability}, the GIT stability. 
It is a quite difficult and interesting problem to 
explicitly understand the stability notion. 

Let us recall that the compact moduli scheme   
of curves $\bar{M_{g}}$ 
is constructed in GIT by permitting ordinary 
double points {\color{black}{(nodes)}} to curves (\cite{DM69}, \cite{KM76}, \cite{Mum77}, \cite{Gie82}), which is sometimes called the 
{\color{black}{Deligne-Mumford compactification}}. 
{\color{black}{We note that 
semistable polarized curves have only nodal singularities. 
}}
 
In this paper, we give its higher dimensional generalization and 
show that the general effect of singularities on stability is determined by the  \textit{discrepancy}, an invariant of singularity which was developed along the 
minimal model program. This is our new point of view. 
{\color{black}{Recall that the discrepancy is defined under the following conditions, 
which ensure that the canonical divisor $K_X$ or the canonical sheaf $\omega_X$ is 
in a tractable class (cf. e.g., \cite{Ale96}). }}

\begin{Def}
{\color{black}{
An algebraic scheme $X$ is said to {\it satisfy $(*)$} when 
the following conditions hold. }}
\begin{enumerate}
{\color{black}{
\item $X$ is equidimensional and reduced. 
\item $X$ satisfies Serre's $S_2$ condition. 
\item Codimension $1$ points of $X$ are Gorenstein. 
\item $K_X$ is $\mathbb{Q}$-Cartier, 
in the sense that $\mathcal{O}_{X}(nK_X):=(\omega_X^{\otimes n})^{\vee \vee}$ 
is an invertible sheaf for some $n\in \mathbb{Z}_{>0}$, where 
$\mathcal{F}^{\vee}:=\mathcal Hom_{\mathcal{O}_X}(\mathcal{F},\mathcal{O}_X)$. }}
\end{enumerate}
\end{Def}

Then, our Main result is the following. 

\begin{Thm}\label{Intro.conj}
\color{black}
{Let $X$ be a projective scheme satisfying $(*)$ and $L$ be an ample 
line bundle on $X$. 
Then, if  $(X,L)$ is $K$-semistable, $X$ has only semi-log-canonical singularities. }
\color{black}{}
\end{Thm}

We will also explain that Theorem \ref{Intro.conj} above should be the best possible 
{\color{black}{as we will see in the statements of Theorem \ref{conv} below.}}  

As already mentioned, the definition of {\it semi-log-canonicity} is based on the 
\textit{discrepancy}. The theory of discrepancy originally 
{\color{black}{stemmed out of a necessity in the way of 
extending minimal models for surfaces by the Italian school to higher dimensions after 
Mori \cite{Mor82}, 
as they should be allowed to have some mild singularities. Indeed, 
it forms a core notion in the minimal model program (the MMP, for short). }}

Along the development of the MMP, 
the semi-log-canonicity {\color{black}{was}} 
first introduced by Koll\'{a}r and Shepherd-Barron 
\cite{KSB88} for {\color{black}{surfaces}} and extended by Alexeev \cite{Ale96} 
{\color{black}{to higher dimensions. Their original purpose was to construct }}
the compactified moduli spaces 
{\color{black}{for varieties of general type}} \textit{not} by GIT theory, 
{\color{black}{but by MMP techniques. }}
For {\color{black}{the case of curves}}, {\color{black}{semi-log-canonical singularities}} are simply smooth points or nodes. Semi-log-canonical surface singularities are classified by 
Koll\'ar-Shepherd-Barron \cite[Theorem (4.24)]{KSB88}. 


Now, let us explain {\color{black}{the other side i.e.,}} the stability notion. 
{\color{black}{While the GIT stability was originally intended to construct moduli spaces 
as mentioned at the beginning, the}} 
\textit{K-}(\textit{semi})\textit{stability} is a version of GIT-stability notion which was  firstly introduced by Tian \cite{Tia97} to describe when 
a Fano manifold has a K\"ahler-Einstein metric. 
{\color{black}{Subsequently}}, Donaldson \cite{Don02} 
{\color{black}{extended the notion to general polarized varieties}} 
with an expectation 
{\color{black}{of correspondence with the}} existence of K\"ahler metrics whose scalar curvature are constant (cscK, for short). 
We follow {\color{black}{Donaldson's formulation}} \cite{Don02} 
in this paper. We note that it is defined {\color{black}{algebro-geometrically}}, 
although the introduction is motivated by differential geometry. 

{\color{black}{
Thus roughly speaking, 
our Main theorem \ref{Intro.conj} bridges in a fresh way, 
these two theories in algebraic geometry, i.e., birational geometry and GIT stability 
(in a broader sense). 
In addition, due to the conjectural correspondence with metrics side, 
one could hope that stability or moduli problems have further connections with differential geometry. 
}}


We should make some comments on Theorem \ref{Intro.conj}. 

{\color{black}{Firstly}}, we {\color{black}{remark on}} 
the Fano case. 
{\color{black}{In this paper, $X$ is said to be  
{\textit{$(*)$-Fano scheme}} if 
$X$ is a projective scheme satisfying $(*)$ and 
$-K_X$ is ample 
(we do \textit{not} \color{black}{a priori}\color{black}{} assume normality of $X$). In this  case, we can prove the following stronger result by slightly different arguments. 

\begin{Thm}\label{Fano}
If $X$ is {\color{black}{a}} {\color{black}{$(*)$-Fano scheme as above}} 
and $(X,\mathcal{O}_{X}(-mK_X))$ is K-semistable 
with $m\in \mathbb{Z}_{>0}$, then $X$ is log terminal. $($In particular, $X$ should be normal$)$. 
\end{Thm}

{\color{black}{Secondly}}, let us comment on other stability notions. 
Recall that Mumford and Gieseker studied 
\textit{asymptotic} (\textit{Chow} and \textit{Hilbert}) \textit{stabilit
{\color{black}{ies}}}, which were the 
original stability notions for polarized varieties (\cite{Gie82}, \cite{Mum77} etc). 
It is well known that these asymptotic (Chow or Hilbert) semistabilit{\color{black}{ies}}  impl{\color{black}{y}} K-semistability (cf.\ \cite[section 3]{RT07}). Furthermore, there are more stability notions introduced recently by Donaldson 
{\color{black}{(\cite{Don10})}}, {\color{black}{called}} $\bar{K}$-stability 
and b-stability. 
{\color{black}{It seems that these two notions are}} expected to be equivalent at least 
for smooth case, and we can see that $\bar{K}$-semistability is also stronger than 
K-semistability. Therefore, we have 

\begin{Cor}\label{otst}
{\rm (i)}
\color{black}
{Let $X$ be a projective scheme satisfying $(*)$ and $L$ be an ample 
line bundle on $X$. }
\color{black}
{Then, if  $(X,L)$ is asymptotically (Chow or Hilbert) semistable, $X$ has only semi-log-canonical singularities. }
\color{black}{}

{\rm (ii)}
\color{black}
{Let $X$ be a projective scheme satisfying $(*)$ and $L$ be an ample 
line bundle on $X$.  }
\color{black}
{Then, if  $(X,L)$ is $\bar{K}$-semistable, $X$ has only semi-log-canonical singularities. }
\color{black}{}

{\rm (iii)}If $X$ is a {\color{black}{$(*)$-Fano scheme}} 
and $(X,\mathcal{O}_{X}(-mK_X))$ 
with $m\in \mathbb{Z}_{>0}$ is asymptotically 
$($Chow or Hilbert$)$ semistable or $\bar{K}$-semistable, then $X$ is log terminal. 
(In particular, $X$ should be normal). 
\end{Cor}

Final but an important remark about Theorem \ref{Intro.conj} 
is that the following converse 
has already been proved for Calabi-Yau case (\cite{Od09}) and canonically polarized case (\cite{Od11}). 
{\color{black}{In this sense, Theorem \ref{Intro.conj} is the best possible as 
mentioned earlier. }}

\begin{Thm}\label{conv}

{\rm (i)}$($\cite{Od09}$)$
A semi-log-canonical polarized variety $(X,L)$ 
with numerically trivial canonical divisor $K_X$ 
is K-semistable. 

{\rm (ii)}$($\cite{Od11}$)$
A semi-log-canonical $($pluri$)$canonically polarized variety $(X,\mathcal{O}_{X}(mK_{X}))$ with $m\in \mathbb{Z}_{>0}$ is K-stable. 
\end{Thm}

\noindent
We make a caution that, on the other hand, 
the singularities {\color{black}{do}} not determine stabilities \textit{in general}, as 
it is well known that there are smooth but 
\color{black}{not semistable}\color{black}{} polarized manifolds. 

{\color{black}{
We also remark that it has been known for a few decades that 
{\color{black}{the}} 
asymptotic (semi)stability version of Theorem \ref{conv} 
does {\it not} {\color{black}{hold}} (cf.\ \cite{She83}, \cite{Od11}). 
The author supposes {\color{black}{that}} 
this phenomenon ought to be a major reason why the relation between 
discrepancy and stability of polarized varieties has been 
{\color{black}{unexpected}} so far. }}


{\color{black}{
The bare structure of the proof of Theorem \ref{Intro.conj} is that, }}
assuming non-semi-log-canonicity of $X$ (i.e., $X$ has ``bad" singularities), 
we construct a {\color{black}{``de-stabilizing" one-parameter subgroup}} 
by making use of a certain 
birational model of $X$ and $X\times\mathbb{A}^1$. 
{\color{black}{On the way of the proof, we define \textit{S-coefficient}, }}
which is an invariant of certain ideals of $X\times \mathbb{A}^{1}$. 
{\color{black}{Very roughly speaking, to those ideals we associate the one parameter subgroups. }}

{\color{black}{The}} birational model of $X$ which we {\color{black}{shall use }}
is the (relative) \textit{semi-log canonical model} 
whose existence has been conjectured in the theory of the 
log minimal model program (LMMP, for short), at least for normal case. 
The existence is recently verified in \cite{OX11}.

In our standpoint, Shah \cite{Sha81} introduced 
{\color{black}{our key invariant 
S-coefficient for isolated singularities 
by an argument based on Eisenbud-Mumford's local stability theory \cite{Mum77}, 
and applied it to give certain list of semistable surface singularities,  
which gave us one of the major inspirations}} for Theorem \ref{Intro.conj}. 

Our paper is organized as follows. 

In the next section, we will review the basic stability notions for polarized varieties 
and some preparatory materials related to the log minimal model program. 
In section $3$, we will formulate an invariant 
of polarized varieties (with an ideal of certain type attached),  
which we call the \textit{S-coefficient}, as a generalization of ``$a_{I}$" in 
\cite{Sha81}. 
Actually, the S-coefficient can be regarded as the leading coefficient of some series of the Donaldson-Futaki invariants, which can be calculated by formula  \ref{DF.formula} proven in \cite{Od09}, \cite{Wan08}. 
After that, we 
{\color{black}{give technical details to the (birational geometric part of) proof of}} 
Theorem \ref{Intro.conj}, \ref{Fano} in the following 
sections. 

\subsection*{Conventions}

Throughout, we work over an algebraically closed field $k$ with 
characteristic $0$.

A polarization means an \textit{ample} invertible sheaf and a 
polarized scheme means {\color{black}
{an algebraic scheme $X$ equipped }}
with an ample invertible sheaf $L$. $(X,L)$ always denotes a polarized scheme, and 
except in subsection \ref{K-stability} 
{\color{black}{and a tiny part of subsection \ref{formula.subsection}, 
it is assumed to satisfy $(*)$ 
as in the statement of Theorem \ref{Intro.conj}. 
}}(For example, an arbitrary 
{\color{black}{reduced}} projective hypersurface, or more generally, 
a (global) {\color{black}{reduced}} 
complete intersection satisfies the condition{\color{black}{s}}. )

$\NN(X)$, $\NLC(X)$, $\NSLC(X)$ and $\NKLT(X)$ denote non-normal locus, non-log-canonical locus, non-semi-log-canonical 
locus, and non-kawamata-log-terminal locus of $X$, respectively. 
$X^{\nu}$ denotes the normalization of a given variety $X$. 

$a({\color{black}{E}};X)$ denotes the discrepancy of 
{\color{black}{a divisor}} ${\color{black}{E}}$ 
{\color{black}{over a}} normal variety $X$ 
and $a({\color{black}{E}};X,D)$ denotes the discrepancy of ${\color{black}{E}}$ over a 
{\color{black}{normal}} pair $(X,D)$ (i.e., a pair of a normal variety $X$ and its Weil divisor $D$ with $\mathbb{Q}$-Cartier $K_{X}+D$). 

\subsection*{Acknowledgements} 

This paper is an expanded version of the author's master thesis, 
\color{black}
{submitted}\color{black}{} in January of 2009, to RIMS, Kyoto university. 
It essentially proved Theorem \ref{Intro.conj} modulo the existence of the 
(relative) semi-log-canonical model. 

The author can not express enough his sincere and deep 
gratitude to his advisor Professor Shigefumi Mori for his 
helpful comments, reading the draft and the revision 
and for his warm encouragements. 
The author also would like to thank Professors 
Shigeru Mukai, Noboru Nakayama, Masayuki Kawakita, Osamu Fujino, 
Julius Ross, Xiaowei Wang, Yuji Sano, Yongnam Lee, Shiro Goto, 
Hisanori Ohashi and Mr. Kento Fujita 
for helpful comments and encouragements for pursuing the program, 
writing {\color{black}{and revising}} the paper. For the revision, 
the anonymous referee considerably helped the author to improve the exposition, 
which he also greatly appreciate. 

My big special thanks go to Professors J\'anos Koll\'ar 
and Chenyang Xu. It is a great pleasure of the author to mention that, 
thanks to them, the technical problem of constructing the relative semi-log-canonical model is settled in \cite{OX11} on which our proof depends. 

The author is 
{\color{black}{partially}} supported by the Grant-in-Aid for Scientific Research (KAKENHI No.\ 21-3748) and the Grant-in-Aid for JSPS fellows.


\section{Preliminaries}

\color{black}
{
In this section, we review the basics for K-stability, 
discrepancy, and log canonical model. 
}\color{black}{}




 

\subsection{K-stability}\label{K-stability}

The K-stability was introduced {\color{black}{first}}, under differential geometric background, by Tian in \cite{Tia97}, and reformulated and extended later by Donaldson \cite{Don02}. 
{\color{black}{
Donaldson's version of K-stability and K-polystability have been slightly amended  recently by 
\cite{LX11}, while the semistability notion remained the same (see also \cite{Od12}, \cite{Stp11}). 
}}
Recall that 
{\color{black}{it is}} the motivation for introducing the K-(semi, poly)stability to 
seek the GIT-counterpart of the existence of special K\"{a}hler metric. 
{\color{black}{
Indeed, according to Professors Gang Tian and Toshiki Mabuchi, 
the ``K" in K-stability stands for the K-energy (Mabuchi energy), 
a functional on the space of K\"ahler metrics whose critical points are 
canonical K\"ahler metrics and at last the ``K" in the K-energy came from 
``K"\"ahler. }}

For the definition of the stability, 
we need the concept of ``test configuration" following Donaldson \cite{Don02}. 
Our notation (and even expression) 
{\color{black}{mostly}} follows \cite{RT07}. 

\begin{Def}

A \textit{test configuration} (resp.\ \textit{semi test configuration}) for a polarized 
\color{black}{complete}\color{black}{} 
scheme $(X,L)$ is a quasi-projective scheme $\mathcal{X}$ with 
an invertible sheaf $\mathcal{M}$ on it with: 
\begin{enumerate}
\item{a $\mathbb{G}_{m}$ action on $(\mathcal{X},\mathcal{M})$}
\item{a proper flat morphism $\alpha\colon \mathcal{X} \rightarrow \mathbb{A}^{1}$}
\end{enumerate}
such that $\alpha$ is $\mathbb{G}_{m}$-equivariant for the usual action on $\mathbb{A}^{1}$: 
\begin{align*}
\mathbb{G}_{m}\times \mathbb{A}^{1}&& \longrightarrow&& \mathbb{A}^{1}\\
                          (t,x)    && \longmapsto    &&    tx,      
\end{align*}
$\mathcal{M}$ is relatively ample (resp.\ relatively semi ample), 
and $(\mathcal{X},\mathcal{M})|_{\alpha^{-1}(\mathbb{A}^{1}\setminus \{0\})}$ is $\mathbb{G}_{m}$-equivariantly isomorphic 
to $(X,L^{\otimes r})\times (\mathbb{A}^{1}\setminus \{0\})$ for some positive integer $r$, called \textit{exponent}, 
with the natural action of $\mathbb{G}_{m}$ on the latter and the trivial action on the former. 

\end{Def}

\begin{Prop}[{\cite[Proposition 3.7]{RT07}}]\label{tc.1-ps}

In the above situation, a one-parameter subgroup of $GL(H^{0}(X,L^{\otimes{r}}))$ is equivalent to the data of a test 
configuration $(\mathcal{X},\mathcal{M})$ of $(X,L)$ with the polarization 
$\mathcal{M}$ very ample $($over $\mathbb{A}^{1}$$)$ and 
{\color{black}{of}} exponent {\color{black}{$r$ for $r\gg 0$. }}

\end{Prop}

\noindent
In fact, let $\lambda \colon \mathbb{G}_{m}\rightarrow \GL(H^{0}(X,L^{\otimes r}))$ 
be a one-parameter subgroup. Then, consider the natural action $\lambda \times \rho$  of $\mathbb{G}_{m}$ on $(\mathbb{P}(H^{0}(X,L^{\otimes r}))\times \mathbb{A}^{1},  \mathcal{O}(1))$ as a polarized variety, 
where  $\rho$ is the multiplication action on $\mathbb{A}^{1}$. 
Then the closure of the orbit $\mathcal{X}:=\overline{((\lambda \times \rho) (\mathbb{G}_{m}))(X \times \{1\})}$ is a test configuration with the natural polarization $\mathcal{O}(1)|_{\mathcal{X}}$ and the restriction of the natural action on 
$(\mathbb{P}(H^{0}(X,L^{\otimes r}))\times \mathbb{A}^{1},  \mathcal{O}(1))$. This is 
called the \textit{DeConcini-Procesi family} of $\lambda$ by Mabuchi. 
The fact that any (very ample) test configuration can be obtained in this way follows from 
\color{black}{the fact that an arbitrary $\mathbb{G}_m$-equivariant vector bundle over $\mathbb{A}^1$ should be equivariantly trivial (cf.\ \cite[Lemma 2]{Don05}). 
}\color{black}{}

Therefore, the test configuration can be regarded as \textit{geometrization} of one-parameter subgroup. 

\color{black}{
Now, let us define the Donaldson-Futaki invariants for test configurations 
whose positivity define K-stability. 
}\color{black}{}
As a preparation, let us note that the \textit{total weight} of an action of $\mathbb{G}_{m}$ on some finite-dimensional vector space will mean 
the sum of all weights in this paper. Here the \textit{weights} mean the exponents of eigenvalues which should be powers of $t$. 
\color{black}{
Take a test configuration $(\mathcal{X},\mathcal{M})$ and 
suppose that the exponent $r$ is $1$. Otherwise, we can similarly 
proceed by considering $(X,L^{\otimes r})$ instead of $(X,L)$. 
}\color{black}{}
We denote the total weight of the induced action on $(\alpha_{*}\mathcal{M}^{\otimes{U}})|_{0}$ {\color{black}{by}} 
$w(U)$ and $\dim X$ as $n$. 
It is a polynomial of $U$ of degree $n+1$. 
\color{black}{
On the other hand, 
}\color{black}{}
we write $P(u):=\dim H^{0}(X,L^{\otimes{u}})$. 
Let us take $rP(r)$-th power \color{black}{
of the action of $\mathbb{G}_{m}$ on 
$\mathcal{M}|_{0}$ and multiply suitable power of $t$ so that 
the action on the vector space $(\alpha_*\mathcal{M}^{\otimes r})|_{\{0\}}$ 
would be in the special linear group ${\rm SL}
((\alpha_*\mathcal{M}^{\otimes r})|_{\{0\}})$.  Then, }\color{black}{} 
the corresponding normalized weight on 
$(\alpha_{*}\mathcal{M}^{\otimes{U}})|_{0}$ 
is $\tilde{w}_{r,Ur}:=w(u)rP(r)-w(r)uP(u)$, where $u:=Ur$. It is a 
polynomial of form 
$\sum_{i=0}^{n+1}e_{i}(r)u^{i}$ of degree $n+1$ in $u$ for $u \gg 0$. 
\color{black}{
Further, the coefficients $e_i(r)$ are again polynomials of degree $n+1$ in $r$ 
for $r \gg 0$ : $e_{i}(r)=\sum_{j=0}^{n+1}e_{i,j}r^{j}$ for $r \gg 0$. 
}\color{black}{}
Since the weight is normalized, 
$e_{n+1,n+1}=0$. $e_{n+1,n}$ is 
called the \textit{Donaldson-Futaki invariant} of the test configuration, which we will denote {\color{black}{by}} 
$\DF(\mathcal{X},\mathcal{M})$. 
\color{black}{
Note that $(n+1)!e_{n+1}(r)r^{n+1}$ has meaning as the Chow weight of $X\subset \mathbb{P}(H^{0}(X,L^{\otimes r}))$ 
with respect to the SL-normalization of the one parameter subgroups associated to 
$(\mathcal{X},\mathcal{M}^{\otimes r})$ via Proposition \ref{tc.1-ps} for $r \gg 0$ 
(cf.\ \cite[Lemma 2.11]{Mum77}). 
}\color{black}{}

For an arbitrary \textit{semi} test configuration $(\mathcal{X},\mathcal{M})$ 
we can define the (normalized) Chow weight or 
the Donaldson-Futaki invariant in completely simlar way 
{\color{black}{from the total weights}} of the induced $\mathbb{G}_m$-action on 
{\color{black}{$(\alpha_*\mathcal{M}^{\otimes U})|_{\{0\}}$ for $U\gg 0$. 
Also note that the homogeneity 
$\DF(\mathcal{X},\mathcal{M}^{\otimes c})=c^{2n}\DF(\mathcal{X},\mathcal{M})$ 
easily follows from the definition. 
}}

Now, we can recall the definition of K-stability as follows. 

\begin{Def}[{cf.\ \cite{Stp11}, \cite{Od12}}]
\color{black}{
A test configuration $(\mathcal{X},\mathcal{L})$ is said to be {\it almost trivial} 
if $\mathcal{X}$ is $\mathbb{G}_m$-equivariantly isomorphic to the product 
test configuration away from a closed subset of codimension at least $2$. 
}\color{black}{}
\end{Def}

\begin{Def}

{\rm (i)} 
A polarized {\color{black}{complete}} scheme $(X,L)$ is \textit{K-stable} 
(resp.\ \textit{K-semistable}) if for any 
{\color{black}{
test configurations {\color{black}{of $(X,L)$}} which are not almost trivial, 
}}
 with exponent $r$, the leading coefficient $e_{n+1,n}$ of $e_{n+1}(r)$ 
(the Donaldson-Futaki invariant) is 
positive (resp.\ non-negative). 

{\rm (ii)} 
{\color{black}{
A polarized {\color{black}{complete}} scheme $(X,L)$ is \textit{K-polystable} 
if it is K-semistable and the Donaldson-Futaki invariant of a test configuration 
$(\mathcal{X},\mathcal{M})$ is $0$ if and only if $\mathcal{X}$ is isomorphic to 
$X\times \mathbb{A}^1$ away from a closed subset of codimension at least $2$. 
}}

\end{Def}

{\color{black}{
Although we only use K-semistability in this paper, we make remarks on 
other notions. }}
We should note that the original ``K-stability" of \cite{Don02} is what is called 
{\color{black}{``K-\textit{poly}stability"}} in \cite{RT07}. 
We follow the convention of \cite{RT07} 
{\color{black}{at this point.}} We further note that it is possible to 
re-define asymptotic stability by the quantities introduced above, associated 
to test configurations, due to Proposition \ref{tc.1-ps}.





\color{black}{
About other stability notions, we only note that stablity notions are related as follows, 
without giving their definitions and proofs. For the details, we refer to 
\cite{RT07} and \cite{Don10}. 
}\color{black}{}

\begin{Claim}\label{compa}

\begin{enumerate}

\item{Asymptotically Chow stable $\Rightarrow$ Asymptotically Hilbert stable 
$\Rightarrow$ Asymptotically Hilbert semistable $\Rightarrow$ Asymptotically Chow semistable 
$\Rightarrow$ K-semistable. }
\item{$\bar{K}$-stable $\Rightarrow$ $\bar{K}$-semistable $\Rightarrow$ K-semistable.  }
\end{enumerate}
\end{Claim}
\noindent
\color{black}{
Hence, among these notions, K-semistability is the weakest notion. 
It is the reason why 
Corollary \ref{otst} should follow from 
Theorem \ref{Intro.conj}, \ref{Fano}. 
}\color{black}{}

\subsection{Singularities via discrepancy}\label{sing.disc}
We will now explain the discrepancy and some 
{\color{black}{classes of mild singularities}}. 
Consult \cite[section 2.3]{KM98} and \cite[Chapter 12]{Koletc92} for the details. 
Let us {\color{black}{first}} treat normal case. Let $(X,D)$ be a normal pair, i.e., a pair of a normal variety $X$ and an effective $\mathbb{Q}$-divisor $D$ such that $K_{X}+D$ is $\mathbb{Q}$-Cartier. Let $\pi:X'\to X$ be a log resolution of $D$, i.e., $\pi$ is a proper birational morphism such that $X'$ is smooth and $\pi^{-1}\Supp(D)\cup E$ 
has a simple normal crossing divisor support, where $E$ is the exceptional divisor of $\pi$. Then, we denote 
\begin{equation*}\label{eq:discrepancy}
	K_{X'}-\pi^*(K_{X}+D)=\sum_i \color{black}{a(E_i;(X,D))}\color{black}{}E_{i}, 
\end{equation*} 
where $a(E_i;(X,D)) \in \mathbb{Q}$ and $E_{i}$ {\color{black}{run}} over the set of divisors of $X'$ supported on the exceptional locus or the support $\Supp(\pi^{-1}_{*}D)$ of 
{\color{black}{$\pi^{-1}_*D$}}, the strict transform of $D$. 
\color{black}{
We sometimes simply write $a(E_i;(X,D))$ as $a(E_i;X)$ if $D=0$,  
and write $a(E_i)$ if the pair in concern is obvious from the context. 
}\color{black}{}

The pair $(X,D)$ is called \textit{log canonical} 
(resp.\ \textit{{\color{black}{kawamata}} log terminal}) 
if and only if $a(E_i;(X,D))\ge -1$ ({\color{black}{resp.\ $a(E_i;(X,D))>-1$}}) for any $E_{i}$. 
These notions are  independent of the choice of the log resolution. 
{\color{black}{
We simply call $X$ is {\it log canonical} (resp.\ {\it log terminal}) 
when $(X,0)$ is log canonical (resp.\ kawamata log terminal). 
}}

The semi-log-canonicity is an extension to non-normal case  of the notion of log-canonicity. 
We introduce those notions without divisors, i.e. in non-log setting, 
{\color{black}{at this stage of argument. }}
(We will need some log versions as well later in section \ref{sec.NN2}, 
where we introduce those definitions.)

Let $X$ be 
a projective variety, which is reduced, equidimensional, 
$\mathbb{Q}$-Gorenstein, 
Gorenstein in codimension $1$ and satisfies Serre condition $S_{2}$ (as we assumed).  Let $\nu \colon X^{\nu}\to X$ be the normalization morphism and 
attach a conductor divisor $\cond(\nu)$ on $X^{\nu}$ which is defined by $K_{X^{\nu}}=\nu^{*}K_{X}+\cond(\nu)$. From the assumption, $(X^{\nu},\cond(\nu))$ is 
a log pair (i.e. $K_{X^{\nu}}+\cond(\nu)$ is $\mathbb{Q}$-Cartier). Then, the 
\textit{semi-log-canonicity} of $X$ are defined 
simply as the log canonicity of the normalized pair, $(X^{\nu},\cond(\nu))$. 

For curve case, the semi-log-canonicity is equivalent to that the curve is nodal (or smooth). For surface case, that class of singularities is also classified by Koll\'ar-Shepherd-Barron \cite{KSB88}. For higher dimensional case, it is well known that 
a semi-log-canonical variety has only normal crossing singularity in codimension $1$, 
so that repeatedly taking general hyperplane section leads to a nodal curve. 

\subsection{Log canonical model}\label{lc.mod.sec}
\color{black}{
To construct ``de-stabilizing" test configurations for non-semi-log-canonical 
polarized varieties, we need a birational model called \textit{(relative) log canonical model}. The definition is as follows. 
}\color{black}{}
\begin{Def}
\color{black}{
Let $(X,D)$ be a normal pair, i.e. $X$ is a normal variety attached 
with a $\mathbb{Q}$-divisor such that $K_X+D$ is $\mathbb{Q}$-Cartier. 
We call that a birational projective morphism $\pi\colon B\to (X,D)$ 
gives  {\it a (relative) log canonical model} of $(X,D)$ 
{\color{black}{(or of $X$ if $D=0$)}} if with the divisor 
$E_{\rm red}$, which denotes the sum of $\pi$-exceptional prime divisors with coefficients $1$, the pair $(B,E_{\rm red})$ satisfies 
}
\begin{enumerate}
\item[(1)] \color{black}{$(B,E_{\rm red})$ is a log canonical pair,}
\item[(2)] \color{black}{$K_B+E_{\rm red}$ is ample over $X$. }
\end{enumerate}
\end{Def}

\noindent
\color{black}{
The existence is established in \cite{OX11}. 
{\color{black}{We used the variable $B$ as we shall use this 
regarding it as a certain blow up of $X$.  
}}
Indeed, this model is a log canonical model of a log resolution with a boundary 
supported on the exceptional set in the sense of log minimal model program. 
}\color{black}{}

\section{The S-coefficients}\label{sec.Scoeff}

In this section, we introduce the concept of {\color{black}{{\it S-coefficients}  
which control asymptotic behaviors for Donaldson-Futaki invariants 
of certain series of test configurations, }}and establish some basic 
properties. 

\subsection{{\color{black}{Review of the formula for Donaldson-Futaki invariants}}}
\label{formula.subsection}

In this subsection, let us recall the formula for Donaldson-Futaki invariants we shall use 
from \cite{Od09}. {\color{black}{
Note that a slightly different version of the formula had been 
also proved independently by Xiaowei Wang in \cite{Wan08}. }}

{\color{black}{Firstly we define a class of ideals, which we shall 
use for our study of stability. }}
Let $(X,L)$ be an $n$-dimensional polarized \color{black}{complete}\color{black}{} variety 
(which is \textit{not} necessarily normal). 

\begin{Def}

A coherent ideal \color{black}{sheaf}\color{black}{} 
$\mathcal{J}$ of $X\times \mathbb{A}^{1}$ is called a \textit{flag ideal} if 
$\mathcal{J}=I_{0}+I_{1}t+\dots+I_{N-1}t^{N-1}+(t^{N})$,  
where $I_{0}\subseteq I_{1}\subseteq \dots I_{N-1} \subseteq \mathcal{O}_{X}$ is the sequence of coherent ideal \color{black}{sheaves}\color{black}{} . 
(It is equivalent to that 
\color{black}{
the corresponding subscheme is 
supported on the central fiber $X\times \{0\}$ and is 
}\color{black}{}
$\mathbb{G}_{m}$-invariant under the natural action 
of $\mathbb{G}_{m}$ on $X\times \mathbb{A}^{1}$.) 
\end{Def}

Let us introduce some notation. 
We set $\bar{\mathcal{L}}:=p_{1}^{*}L$ on $X\times \mathbb{P}^{1}$ and 
its restriction $\mathcal{L}:=p_{1}^{*}L|_{(X\times \mathbb{A}^{1})}$, 
where $p_{i}$ is the $i$-th projection morphism from $X \times \mathbb{A}^{1}$ or $X \times \mathbb{P}^{1}$. 
Let us write the blow up $\bar{\mathcal{B}}(:=Bl_{\mathcal{J}}(X\times \mathbb{P}^{1}))\rightarrow X\times \mathbb{P}^{1}$ 
\color{black}{
or its restriction to $\mathcal{B}(:=Bl_{\mathcal{J}}(X\times \mathbb{A}^{1}))
\rightarrow X\times \mathbb{A}^1$ by $\Pi$. 
Its natural exceptional divisor will be written as $E$
}\color{black}{}, 
i.e.\ $\mathcal{O}_{\mathcal{B}}(-E')=\Pi^{-1}\mathcal{J}$.  
{\color{black}{(We shall use the symbol (prime) $'$ for denoting exceptional divisors to  indicate it they are exceptional divisors of $(n+1)$-dimensional variety $X\times \mathbb{A}^1$, not of $X$.) }}

Let us assume $r$ is sufficiently large so that 
$\color{black}{(\Pi^*\mathcal{L}^{\otimes r})}\color{black}{}(-E')$ is (relatively) semi-ample (over $\mathbb{A}^{1}$). 
Consider the Donaldson-Futaki invariant of the (semi) test configuration 
$(\mathcal{B}, (\Pi^*\mathcal{L})^{\otimes r}(-E'))$. 
Let us recall our formula for that. 

\begin{Thm}[{\cite[Theorem 3.2]{Od09}}]\label{DF.formula} 

Let $(X,L)$ and $\mathcal{B}$, $\mathcal{J}$ be as above. And we assume that exponent $r=1$. 
$($It is just to make the formula easier. For general $r$, put $L^{\otimes r}$ and $\bar{\mathcal{L}}^{\otimes r}$ 
to the place of $L$ and $\bar{\mathcal{L}}$. $)$ 
\color{black}{
Furthermore, we assume that 
$X$ is equidimensional, reduced, satisfying $S_2$ condition, whose codimension 1 points are Gorenstein and having $\mathbb{Q}$-Cartier canonical divisor $K_X$ 
}\color{black}{}
and $\mathcal{B}$ is Gorenstein in codimension $1$. 
Then the corresponding Donaldson-Futaki invariant $\DF((Bl_{\mathcal{J}}(X\times \mathbb{A}^{1}), \mathcal{L}(-E')))$  is 

\begin{equation*}
\dfrac{1}{2(n!)((n+1)!)}\bigl\{-n(L^{n-1}.K_{X})((\Pi^*\bar{\mathcal{L}})(-E'))^{n+1}
\end{equation*}
\begin{equation*}
+(n+1)(L^{n})
(((\Pi^*\bar{\mathcal{L}})(-E'))^{n}.\Pi^{*}(p_{1}^{*}K_{X}))
\end{equation*}
\begin{equation*}
+(n+1)(L^{n})(((\Pi^*\bar{\mathcal{L}})(-E'))^{n}.K_{\bar{\mathcal{B}}/X\times \mathbb{P}^{1}})\bigr\}. 
\end{equation*}
\color{black}{
In the above, the intersection numbers $(L^{n-1}.K_X)$ and $(L^n)$ 
are taken on $X$. On the other hand, $K_{\bar{\mathcal{B}}/X\times \mathbb{P}^1}
:=K_{\bar{\mathcal{B}}}-\Pi^*K_{X\times \mathbb{P}^1}$ is an exceptional divisor 
on $\bar{\mathcal{B}}$ and thus 
$(((\Pi^*\bar{\mathcal{L}})(-E'))^{n}.\Pi^{*}(p_{1}^{*}K_{X}))$ 
and $(((\Pi^*\bar{\mathcal{L}})(-E'))^{n}.K_{\bar{\mathcal{B}}/X\times \mathbb{P}^{1}})$ 
are intersection numbers taken on $\bar{\mathcal{B}}$. 
}\color{black}{}
\end{Thm}
We call the sum of first two terms {\color{black}{the}} \textit{canonical divisor part} 
since they involve intersection numbers with $K_{X}$ or its pullback, 
and the last term will be called {\color{black}{the}} \textit{discrepancy term} since it reflects discrepancies over $X$. 
{\color{black}{We remark that although not all semi test configurations 
are of the form $(\mathcal{B},(\Pi^*\mathcal{L})^{\otimes r}(-E'))$, 
it is sufficient for K-(semi)stability to check the Donaldson-Futaki invariants 
of the special semi test configurations (\cite{Od09}). }}


\subsection{{\color{black}{S-coefficient as a leading coefficient of Donaldson-Futaki invariants}}}

We define the S-coefficient, the key invariant as follows. 

\begin{Def}\label{S-coeff.def}

Let us fix $(X,L)$ in Theorem \ref{DF.formula} above and fix a flag ideal $\mathcal{J}$. 
{\color{black}{
Suppose that $\bar{\mathcal{B}}$ is Gorenstein in codimension $1$ so that 
the canonical divisor class $K_{\bar{\mathcal{B}}}$ is well defined. 
}}
Then, the \textit{S-coefficient} for that flag ideal $\mathcal{J}$ 
is defined as an intersection number $(\mathcal{L}^{s}.(-E')^{n-s}.K_{\bar{\mathcal{B}}/(X\times\mathbb{P}^{1})})$ taken on $\bar{\mathcal{B}}$ 
and we denote it {\color{black}{by}} $S_{(X,L)}(\mathcal{J})$, where $s$ denotes the dimension of $\Supp(\mathcal{O}_{X\times 
\mathbb{A}^{1}}/\mathcal{J})$.   
We note that \color{black}{
homogeneity
}\color{black}{} 
$S_{(X,L^{\lambda_{1}})}(\mathcal{J}^{\lambda_{2}})=\lambda_{1}^{s}\lambda_{2}^{n-s}S_{(X,L)}(\mathcal{J})$ 
follows from the definition. 

\end{Def}

\noindent
The main motivation for above definition is the following meaning 
of S-coefficient, as leading coefficient of Donaldson-Futaki invariants. 

\begin{Prop}\label{S-coeff.fund.thm}

Let $(X,L)$ and $\mathcal{J}$ be as above. 
\color{black}{
Then, the following holds. }

\color{black}{
{\rm (i)} The sequence of Donaldson-Futaki invariants  
$\DF(Bl_{\mathcal{J}}(X\times\mathbb{A}^{1}),\mathcal{L}^{\otimes{r}}(-E'))$ for $r \gg 0$,  forms a polynomial. }

\color{black}{
{\rm (ii)}
Its coefficient of $r^{d}$ is $0$ for $d>n+s$ and equals to 
$$
\dfrac{\binom{n}{s} (L^{n})}{2(n!)^{2}}S_{(X,L)}(\mathcal{J})
$$ for $d=n+s$. }

\color{black}{
Hence, if $S_{(X,L)}(\mathcal{J})<0$ for some flag ideal $\mathcal{J}$,  
then $(X,L)$ is \color{black}{not K-semistable}
}\color{black}{}. 

\end{Prop} 

{\color{black}{To prove Proposition \ref{S-coeff.fund.thm} above and analyze the positivity of the S-coefficients later, we shall use the following general properties of intersection numbers. As it follows from a standard arguments, we omit the proof. However, we give statements here for the readers' convenience as it shall be a key for our estimation. }}

\begin{Lem}\label{int.number}

{\color{black}{
Let $\mathcal{X}$ be an arbitrary $n+1$-dimensional equidimensional complete scheme, and $\pi\colon 
\bar{\mathcal{B}}\rightarrow \mathcal{X}$ a surjective, generically finite morphism. 
Then 
}}

{\color{black}{
{\rm (i)}
\begin{equation*}
(\pi^{*}D_{1}.\ldots.\pi^{*}D_{s}.E'_{1}.\ldots.E'_{n+1-s})=0
\end{equation*} 
for arbitrary Cartier divisors $D_{1},\ldots,D_{s}$ on $\mathcal{X}$, and arbitrary Cartier divisors $E'_{1},\ldots,E'_{n+1-s}$ 
with $\dim(\pi(\cap{\Supp(E'_{l})}))<s$. 
}}

{\color{black}{
{\rm (ii)}
\begin{equation*}
(\pi^{*}D_{1}.\ldots.\pi^{*}D_{s}.E'_{1}.\ldots.E'_{n+1-s})>0
\end{equation*} 
for arbitrary ample Cartier divisors $D_{1},\ldots,D_{s}$ on $\mathcal{X}$, 
arbitrary ample Cartier divisors $E'_{1},\ldots,E'_{n-s}$ on $\bar{\mathcal{B}}$ and 
an arbitrary effective Weil divisor $E'_{n+1-s}$ on $\bar{\mathcal{B}}$ 
with $\dim(\pi(E'_{n+1-s}))=s$. 
}}
\end{Lem}




\begin{proof}[Proof of {\em{Proposition \ref{S-coeff.fund.thm}}}]
{\color{black}{
Replacing $L$ by $L^{\otimes r}$ and 
$\mathcal{L}$ by $\bar{\mathcal{L}}^{\otimes r}$ for the formula \ref{DF.formula}, 
we have the formula of 
$\DF(Bl_{\mathcal{J}}(X\times\mathbb{A}^{1}),\mathcal{L}^{\otimes{r}}(-E'))$. 
From that, Proposition \ref{S-coeff.fund.thm} {\rm (i)} easily follows. 
}}

{\color{black}{
Further, Lemma \ref{int.number} {\rm (i)} 
applied to $\pi=\Pi\colon Bl_\mathcal{J}(X\times \mathbb{P}^1)
\to X\times \mathbb{P}^1$ by taking $D_i:=H\times \mathbb{P}^1$ where $H \in |L^{\otimes m}|$ ($m\in \mathbb{Z}_{>0}$), $E'_i=E'$ for $i\leq n-s$ and 
$E'_{n+1-s}=K_{\bar{\mathcal{B}}/(X\times \mathbb{P}^1)}$, 
implies Proposition \ref{S-coeff.fund.thm} {\rm (ii)} straightforward. 
}}
\end{proof}

\subsection{{\color{black}{S-coefficients and discrepancy}}}

{\color{black}{
In this subsection, 
we shall show a criterion on positivity of S-coefficients, 
which gives a relation with discrepancy. 
}}

Let us assume, from now on, that $X$ is an equidimensional reduced projective variety, satisfies $S_2$ condition and whose codimension $1$ points are Gorenstein. Thus we can define the Weil divisor class $K_X$ which we assume to be $\mathbb{Q}$-Cartier. 
If all codimension $1$ points of $\mathcal{B}$ are Gorenstein, we set  $K_{\bar{\mathcal{B}}/(X\times \mathbb{P}^{1})}:=K_{\bar{\mathcal{B}}}-\Pi^*(K_X\times \mathbb{P}^1)=\sum{\color{black}{a(E'_i)}\color{black}{}}E'_{i}$. 

\begin{Prop}\label{disc.S}
{\color{black}{
Let $X$ be as above and $L$ be an ample line bundle on $X$. 
Moreover, assume that there is a flag ideal $\mathcal{J}$ whose blow up $\mathcal{B}$ 
is Gorenstein in codimension $1$ as noted above. Furthermore, assume that the discrepancies $a(E'_i)$ satisfy the following. $a(E'_i)\leq 0$ for all the exceptional prime divisors $E'_{i}$ on $\mathcal{B}$ which 
dominate $s$ $($maximal$)$-dimensional components 
of $\Supp(\mathcal{O}/\mathcal{J})$ 
and moreover there exists at least one such $i$ with $a(E'_i)<0$. 
Then, we have $S_{(X,L)}(\mathcal{J})<0$. 
}}
\end{Prop}

{\color{black}{
Hence, by combining with Proposition \ref{S-coeff.fund.thm}, 
we have the following criterion for when a polarized variety can be not K-semistable. 
}}

\begin{Cor}\label{stadis.Cor}
{\color{black}{
Let $X$ be as above and 
assume that there is a flag ideal $\mathcal{J}$ whose blow up $\mathcal{B}$ is 
Gorenstein in codimension $1$ and the discrepancies $a(E_i)$ satisfy the following. $a(E'_i)\leq 0$ for all the exceptional prime divisors $E'_{i}$ on $\mathcal{B}$ which 
dominate $s$ $($maximal$)$-dimensional components 
of $\Supp(\mathcal{O}/\mathcal{J})$ 
and moreover there exists at least one such $i$ with $a(E'_i)<0$. 
Then $(X,L)$ is 
not K-semistable 
for an arbitrary polarization $L$. 
}}
\end{Cor}

\begin{proof}[Proof of Proposition \ref{disc.S}] 
{\color{black}{
We have 

\begin{eqnarray}
	\nonumber
		S_{(X,L)}(\mathcal{J})
	&:=&
		(\mathcal{L}^{s}.(-E')^{n-s}.K_{\bar{\mathcal{B}}/(X\times\mathbb{P}^{1})})
	\\
	\label{sc}
	&=&
		(\mathcal{L}^{s}.(\bar{\mathcal{L}}^{\otimes r}-E')^{n-s}.K_{\bar{\mathcal{B}}/(X\times\mathbb{P}^{1})}). 
\end{eqnarray}
Indeed, the equality (\ref{sc}) follows from Lemma \ref{int.number} {\rm (i)} 
applied to $\pi=\Pi\colon Bl_\mathcal{J}(X\times \mathbb{P}^1)
\to X\times \mathbb{P}^1$ by taking $D_i:=H\times \mathbb{P}^1$ 
for $i\leq s+1$ where $H \in |L^{\otimes m}|$ ($m\in \mathbb{Z}_{>0}$) and 
for $i>s+1$, $E'_i=\Pi^*(H\times \mathbb{P}^1)$ or $E'_i=E'$ or $E'_i=K_{\bar{\mathcal{B}}/(X\times \mathbb{P}^1)}$. 

Moreover, the last term 
$(\mathcal{L}^{s}.(\bar{\mathcal{L}}^{\otimes r}-E')^{n-s}.K_{\bar{\mathcal{B}}/(X\times\mathbb{P}^{1})})$ 
is positive due to Lemma \ref{int.number} {\rm (ii)} 
applied to $\Pi\colon Bl_\mathcal{J}(X\times \mathbb{P}^1)
\to X\times \mathbb{P}^1$ again by taking, this time, $D_i:=H\times \mathbb{P}^1$ 
for $i\leq s$, $E'_i (i\leq n-s)$ to be an ample compactification of an ample divisor which belongs to $|(\Pi^*\mathcal{L})^{\otimes r}(-E')|$ on $\mathcal{B}$ to  $\bar{\mathcal{B}}$ 
with $r>1$, and $E'_{n+1-s}:=K_{\bar{\mathcal{B}}/(X\times\mathbb{P}^1)}$. 
}}
\end{proof}

\section{Normal case}\label{norm}


{\color{black}{
As an application of the theory of S-coeffiecients prepared in the previous section, }}
we partially prove Theorem \ref{Intro.conj} for normal case in this section. 
More precisely, let $X$ be a normal variety of pure dimension $n$, 
\color{black}{
having $\mathbb{Q}$-Cartier canonical divisor 
}\color{black}{}
in this section.

\begin{proof}[Proof of {\em{Theorem \ref{Intro.conj}}} for normal $X$]

Thanks to Corollary \ref{stadis.Cor}, it is sufficient to construct 
a flag ideal $\mathcal{J}$ satisfying the following \color{black}{property}\color{black}{}. 

\begin{Property}\label{dim+1.cond}
The blow up $\mathcal{B}$ of $\mathcal{J}$ is normal. 
Furthermore, if we 
let $K_{\mathcal{B}/X\times \mathbb{A}^{1}}=\sum a(E'_i)E'_{i}$, 
then we have $a(E'_i)< 0$ for the discrepancy for an arbitrary 
\color{black}{
{\color{black}{$\Pi$-}}exceptional divisor
}\color{black}{}
$E'_{i}$. 
\end{Property}

We will construct such $\mathcal{J}$ in the following $2$ steps. 
Without loss of generality, we can assume that $X$ is irreducible. 

\begin{Step}
Firstly, we construct a coherent ideal 
\color{black}{sheaf}\color{black}{} 
$I$ of $X$, satisfying the following \color{black}{property}\color{black}{}. 
We denote the blow up of $X$ along $I$ {\color{black}{by}} $\pi\colon B=Bl_{I}(X) \rightarrow X$. 

\begin{Property}\label{dim+0.cond}
The blow up $B$ is normal. Furthermore, if let $s$ be $\dim(\Supp(\mathcal{O}_{X}/I))$, 
then, we have $a(E_{i}; X)< -1$ for the discrepancy for an arbitrary 
{\color{black}{$\pi$-}}exceptional divisor \color{black}{$E_{i}$}\color{black}{}. 
\end{Property} 

\color{black}
{We construct such $I$, using the (relative) log canonical model 
{\color{black}{ (cf.\ subsection \ref{lc.mod.sec})}}
as follows. 
Suppose $\pi \colon B\rightarrow X$ is the (relative) log canonical model of $X$, 
which exists {\color{black}{due to \cite[Theorem1.1]{OX11}}}. 
}
\color{black}{
Then, we take the coherent ideal sheaf 
 $I:=(\pi)_{*}\mathcal{O}_{B}\bigl(m(K_{B/X}+E_{\rm red})      
\bigr)$ for sufficiently divisible $m\in \mathbb{Z}_{>0}$ and the total exceptional divisor $E_{\rm red}$,  then $Bl_{I}(X)\simeq B$. Therefore, this $I$ satisfies Property  \ref{dim+0.cond}. 
}\color{black}{}

\end{Step}
\begin{Step}\label{Step2}
{\color{black}{
Next step starts with taking $I$}} constructed in the previous step. 
Using this, we will construct the flag ideal $\mathcal{J}$ satisfying 
{\color{black}{Property}} \ref{dim+1.cond} as follows. 
From the construction, 
{\color{black}{we have}} $\dim(\Supp(\mathcal{O}_{X}/I))\leq \dim(X)-2$. 
Let us take sufficiently divisible positive integers $m, N$ 
and let us define $\mathcal{J}:=\overline{(I+(t^m))^{N}}$ where the overline denotes the integral closure of the coherent ideal. 
Since it is an invariant ideal with respect to the natural $\mathbb{G}_{m}$ action 
on $X\times \mathbb{A}^{1}$, $\mathcal{J}$ is a flag ideal as well. 
We note that $\mathcal{C}:=Bl_{I+(t)}(X\times \mathbb{A}^1)$ is 
{\it the deformation to the normal cone} (cf.\ \cite{Ful84}, \cite{RT07}) 
but simply taking it is {\it not} sufficient for our purpose in general. 
Geometrically speaking, to take $I+(t^m)$ as above, instead of the simplest $I+(t)$,  
corresponds to take base change of $\mathcal{C}$ 
by $m$-th roots of $t$ (i.e. $s \mapsto t:=s^m$) and 
$\mathcal{B}:=Bl_{\mathcal{J}}(X\times \mathbb{A}^{1})$  
is the normalization of the base change (cf.\ \cite{Vas05}). 

Let us think of the more detailed 
geometric structure of the deformation to the normal cone 
$\mathcal{C}$ and its modification $\mathcal{B}$. 
We know that its central fiber consists of two parts: 
the strict transform of $X\times \{0\}$ canonically isomorphic to $B=Bl_{I}(X)$ (we will identify them from now on), 
and the exceptional divisors $F'_i$ which intersect as $F'_i\cap B=E_i$ 
whose generic points $\eta_i$ are regular. 
Indeed, \'etale locally we can write $t=xy^{c_i}$ 
with \'etale local coordinates (i.e., regular parameters) $x, y$ such that $(x=0)$ corresponds to $B$ and $(y=0)$ corresponds to $F'_i$. 

Based on the above facts, we obtain an \'etale local description of $\mathcal{B}\to \mathcal{C}$ explicitly around the generic point $\eta_i$ of $F'_i\cap B$ as follows. 
We can take an \'etale local coordinate system $(u,y,z_1,\cdots,z_{n-2},s)$ of $\mathcal{B}$ around $\eta_i$, and that of $\mathcal{C}$: $(x,y,z_1,\cdots,z_{n-2},t)$ around the fiber of $\eta_i$ 
which are connected by the following equations. 
Here, $t$ denotes the original coordinate of $\mathcal{C}$ corresponding to the  $\mathbb{A}^1$ direction. 
$$x=u^{c_i}, t=s^m. $$ 
\noindent
We denote the preimage of $F'_i$ by $E'_i$, which is irreducible. 
Then, from the above local description, it directly follows that: 
\begin{equation}\label{Cl.4.3}
a(E'_{i}; X\times \mathbb{A}^{1})=b_{i}\bigl(a(E_{i}; X)+1\bigr), 
\end{equation}
\noindent
where each 
$b_{i}:=\frac{m}{c_{i}}$ is a positive integer as $m$ is sufficiently divisible. 
Therefore, $a(E'_{i}; X\times \mathbb{A}^{1})<0$ follows from Property \ref{dim+0.cond} 
in the previous step of construction. This completed the proof of 
Theorem \ref{Intro.conj} for normal varieties' case. 

\end{Step}

\end{proof}

\section{Non-normal case}\label{sec.NN2}

To {\color{black}{give a}} proof of the main theorem \ref{Intro.conj} in full generality, 
we introduce a non-normal generalization of the 
{\color{black}{(relative)}} log-canonical model, 
{\color{black}{which we used in the previous section for normal case. }}
\color{black}{
A reduced equidimensional variety $X$ is called {\it demi-normal} if $X$ is $S_2$, whose codimension 1 points are regular or ordinary nodes. 
}\color{black}{}

\begin{Def}\label{nonnormalmodel}
{\color{black}{
Let $X$ be a demi-normal projective variety. We call a biratonal projective morphism
 $\pi\colon B\to X$ a {\it {\color{black}{(relative)}} semi-log-canonical model}
if $\pi$ is isomorphic over open locus of $X$ with complement's codimension greater than $1$, and satisfies the following two conditions. 
Here, $E_{\rm red}$ denotes the sum of $\pi$-exceptional prime divisors 
with coefficients $1$. }}
\begin{enumerate}
\item[(1)] {\color{black}{$(B,E_{\rm red})$ is a semi-log-canonical pair. }}
\item[(2)] {\color{black}{$K_B+E_{\rm red}$ is ample over $X$. }}
\end{enumerate}
\end{Def}

\noindent
The existence of such models for any $X$ is {\color{black}{again}} proven in 
\cite[Corollary 1.3]{OX11}. 
Given this birational model, the proof of Theorem \ref{Intro.conj} 
{\color{black}{below}} is similar to the case where $X$ is normal. 

\begin{proof}[Proof of Theorem \ref{Intro.conj}]

Take the (relative) semi-log-canonical model $\pi \colon B\rightarrow X$ 
of $X$, which exists {\color{black}{due to \cite[Corollary 1.3]{OX11}}}. 
Here, we note that all the generic points of {\color{black}{$\pi$-exceptional divisors 
are regular, by the definition of the model. }}
Then, if we apply the negativity lemma \cite[Lemma (3.39)]{KM98} to 
these normalizations, we have ${\color{black}{a_{i}}}<-1$ for any $i$, where $K_{B/X}=\sum 
{\color{black}{a_{i}}}
E_{i}$. 
Therefore, if we take $I:=\pi_{*}(\omega_{B/X}^{[l]}(lE))^{**}$ 
\textcolor{black}{with sufficiently divisible positive integer $l$}, 
\color{black}{
where $E_{\rm red}:=\sum E_i$ denotes the total exceptional divisor of $\pi$, }
\color{black}{}
it would be a \color{black}{coherent ideal sheaf}\color{black}{} 
by Serre's $S_{2}$-condition of $X$. Further, it satisfies 
$Bl_{I}(X)\cong B$ by the relative ampleness of $K_{B/X}+E_{\rm red}$. 

Let us consider a flag ideal $\mathcal{J}'=I+(t^{m})$ on $X\times \mathbb{A}^{1}$ 
for sufficiently divisible positive integer $m$, \textcolor{black}{its} blow up 
$\mathcal{C}=Bl_{\mathcal{J}'}(X\times \mathbb{A}^{1})$ 
and its normalization $\mathcal{C}^{\nu} \rightarrow \mathcal{C}$. 
We denote $\Pi \colon \mathcal{C}^{\nu} \to 
X^{\nu}\times \mathbb{A}^{1}$ the associated morphism. 
We can prove $K_{\mathcal{C}^{\nu}}-\Pi^{*}
(K_{X^{\nu}}\times \mathbb{A}^{1}+\cond(\nu)\times \mathbb{A}^{1})=\sum {\color{black}{a'_{i}}}H'_{i}$ 
with ${\color{black}{a'_{i}}}=b_{i}({\color{black}{a_{i}}}+1)<0$ where $b_{i}$ are some positive integers {\color{black}{for each exceptional divisor $H'_i$}} 
{\color{black}{and $\cond(\nu)$ is the conductor divisor of the normalization. 
The proof is in a completly similar manner as in the previous section, 
by taking $(X^{\nu},\cond(\nu))$ instead of $X$ with the normality assumption. }}

We use the \textit{partial normalization $\mathcal{B}$ of 
$\mathcal{C}$} {\color{black}{which was}} introduced in the proof of \cite[Proposition 3.8]{Od09}. 
The definition is $\mathcal{B}:=\Shpec_{\mathcal{O}_{\mathcal{C}}}(i_{*}\mathcal{O}_{X\times (\mathbb{A}\setminus\{0\})}\cap \mathcal{O}_{\mathcal{C}^{\nu}})$, 
where $i \colon X\times (\mathbb{A}^{1}\setminus \{0\}) \hookrightarrow X\times \mathbb{A}^{1}$ is the open immersion.  
Let $f$ be the associated morphism from 
{\color{black}{$\mathcal{C}^{\nu}$ to $\mathcal{B}$}}. 
{\color{black}{Completely similarly}} as we argued in the {\color{black}{
former half of Step $2$ of the previous section}}, we can take a flag ideal $\mathcal{J}$ 
whose blow up is $\Pi \colon \mathcal{B}=Bl_{\mathcal{J}}(X\times \mathbb{A}^{1})
\rightarrow X\times \mathbb{A}^{1}$. Let us recall the following lemma. 

\begin{Lem}[{\cite[Lemma 3.9]{Od09}}]
The morphism {\color{black}{$f \colon \mathcal{C}^{\nu}\to \mathcal{B}$}} is an isomorphism over an open neighborhood of the generic points of the central fiber. 
\end{Lem}

\noindent
Thus, similarly as in the comparison of discrepancies (\ref{Cl.4.3}), we have 
$K_{\mathcal{B}/X\times \mathbb{A}^{1}}=\sum {\color{black}{a'_{i}}}E'_{i}$ 
with ${\color{black}{a'_{i}}}=b_{i}({\color{black}{a_{i}}}+1)<0$ where $b_{i}$ are the positive integers introduced above, 
and {\color{black}{$E'_i$ is the strict transform of $H'_i$. }}

Therefore, we {\color{black}{complete}} the proof 
{\color{black}{Theorem \ref{Intro.conj}}} thanks to 
{\color{black}{Corollary \ref{stadis.Cor}}}. 

\end{proof}

\section{Fano case}

{\color{black}{
In this section, we prove Theorem \ref{Fano}. 
}}

\begin{proof}[proof of Theorem \ref{Fano}]

In this section we do not use the notion of S-coefficients but the proof is done 
by analyzing the formula 
{\color{black}{for}} the Donaldson-Futaki invariants \ref{DF.formula} more directly. 

Let us take a flag ideal $\mathcal{J}:=\overline{(I+(t))^{N}}$, where $I\subset \mathcal{O}_X$ 
corresponds to the reduced subscheme supported on $\NN(X)$, the non-normal locus of $X$, and $N$ is a sufficiently divisible positive integer. 
{\color{black}{Here, we do not take parameter $m$ into account. }}
We note that $\NN(X)$ is 
\color{black}{
purely codimension $1$ 
{\color{black}{ in $X$}} and  
{\color{black}{their generic points are ordinary double points. }}  
In the case of curves (i.e., $\dim(X)=1$), this means $X$ should be nodal. }
\color{black}{} 
Consider the (semi) test configuration of the blow up type $(\mathcal{B},(\Pi^*\mathcal{L})^{\otimes r}(-E'))$ for that flag ideal $\mathcal{J}$ as we did. 
Then, the S-coefficients become $0$ and the leading coefficients of 
$\DF(\mathcal{B},(\Pi^*\mathcal{L})^{\otimes r}(-E'))$ with respect to the variable $r$ has the same signature as $((\Pi^*\bar{\mathcal{L}})^{n-1}.E'^2)$. This can be shown to be negative by cutting $X$ for $s:=\dim(\Supp(\mathcal{O}_{X\times \mathbb{A}^{1}}/\mathcal{J}))=n-1$ times by general 
{\color{black}{hypersurface sections in $|L^{\otimes m}|$}} for $m\gg 0$, 
reducing to the $n=1$ case. 

Thus, we can assume that $X$ is normal. Let us assume that $X$ is 
{\color{black}{log canonical but not log terminal (i.e., strictly log canonical)}} 
and derives a contradiction. 
{\color{black}{In the sense of log minimal model program,}} a log resolution 
with kawamata-log-terminal boundary $(\tilde{X},(1-\epsilon)E_{\rm red})$ with 
$0<\epsilon \ll 1$ should have a log canonical model $B$ over $X$, 
by 
{\color{black}{\cite[Theorem 1.2]{BCHM09}. 
Note that $B$ should be log terminal and so the morphism $B\rightarrow X$ 
is not isomorphism, which is again a blow up of certain coherent ideal sheaf $I$. 
We further remark that the model of subsection \ref{lc.mod.sec} corresponds to the 
$\epsilon=0$ case. Similarly}} as in section \ref{norm}, 
we construct a flag ideal $\mathcal{J}:=\overline{(I+(t^m))^N}$ 
{\color{black}{where $m, N$ are sufficiently divisible positive integers}}, 
and its blow up $\mathcal{B}:=Bl_{\mathcal{J}}(X\times \mathbb{A}^1)$. 
Then, $K_{\mathcal{B}/X\times \mathbb{A}^{1}}=0$ so that the discrepancy 
term vanishes. 

On the other hand, 
{\color{black}{as $s<n$,}} the canonical divisor part of the formula 
\ref{DF.formula} is negative by \cite[proof of Theorem 2.13]{Od11}. 
{\color{black}{Hence,}} $(X,L)$ should 
be {\color{black}{not K-semistable}}. 
This completes the proof of Theorem \ref{Fano}. 

\end{proof}

\end{document}